\documentclass[1p]{elsarticle}
\usepackage{amssymb}
\usepackage{amsfonts}

\newtheorem{theorem}{Theorem}
\newtheorem{conjecture}[theorem]{Conjecture}
\newtheorem{corollary}[theorem]{Corollary}

\newtheorem{observation}[theorem]{Observation}

\newproof{pf}{Proof}

\begin{document}

\title{Asymptotic confirmation of the Faudree-Lehel Conjecture on irregularity strength for all but extreme degrees}

\author[agh]{Jakub Przyby{\l}o\fnref{MNiSW}}

\fntext[MNiSW]{This work was partially supported by the Faculty of Applied Mathematics AGH UST statutory tasks within subsidy of Ministry of Science and Higher Education.}

\address[agh]{AGH University of Science and Technology, Faculty of Applied Mathematics, al. A. Mickiewicza 30, 30-059 Krakow, Poland}

\begin{abstract}
The irregularity strength of a graph $G$, $s(G)$, is the least $k$ admitting a $\{1,2,\ldots,k\}$-weighting of the edges of $G$ assuring distinct weighted degrees of all vertices, or equivalently the least possible maximal edge multiplicity in an irregular multigraph obtained of $G$ via multiplying some of its edges. The most well-known open problem concerning this graph invariant is the conjecture posed in 1987 by Faudree and Lehel that there exists a constant $C$ such that $s(G)\leq \frac{n}{d}+C$ for each $d$-regular graph $G$ with $n$ vertices and $d\geq 2$ (while a straightforward counting argument yields $s(G)\geq \frac{n+d-1}{d}$). The best known results towards this imply that $s(G)\leq 6\lceil\frac{n}{d}\rceil$ for every $d$-regular graph $G$ with $n$ vertices and $d\geq 2$, while $s(G)\leq (4+o(1))\frac{n}{d}+4$ if $d\geq n^{0.5}\ln n$. 

We show that the conjecture of Faudree and Lehel holds asymptotically in the cases when $d$ is neither very small nor very close to $n$. We in particular prove that for large enough $n$ and $d\in [\ln^8n,\frac{n}{\ln^3 n}]$, $s(G)\leq \frac{n}{d}(1+\frac{8}{\ln n})$, and thereby we show that $s(G) = \frac{n}{d}(1+o(1))$ then. We  moreover prove the latter to hold already when $d\in [\ln^{1+\varepsilon}n,\frac{n}{\ln^\varepsilon n}]$ where $\varepsilon$ is an arbitrary positive constant.

\end{abstract}


\maketitle

\section{Introduction}

One of the most basic arguments in graph theory is the pigeonhole principle based observation that the vertices of a simple graph cannot all have pairwise distinct degrees, in other words no graph is in this sense \emph{irregular} except the trivial one-vertex case. This exclusion considered, Chartrand, Erd\H{o}s and Oellermann proposed and investigated in~\cite{ChartrandErdosOellermann} possible alternative definitions of irregularity in graphs. 
In reference to this research, Chartrand et al. gave rise in~\cite{Chartrand} to a related concept, designed to measure in some sense the level of irregularity of a given graph, exploiting for this aim the fact that irregular \emph{multigraphs} are quite common and relatively easily constructed. Namely, they defined the \emph{irregularity strength} of a graph $G=(V,E)$ as the least positive integer $k$ such that one may obtain an irregular multigraph (i.e. a multigraph with pairwise different vertex degrees) by blowing each edge of $G$ up to at most $k$ parallel edges. This graph invariant is usually denoted by $s(G)$ and is undefined for graphs with an isolated edge or two isolated vertices.
An equivalent and often more convenient definition of $s(G)$ relates with so-called \emph{weighted degrees}, defined for a given edge weighting function $\omega: E\to \mathbb{R}$ and a vertex $v\in V$ as:
$$\sigma_{\omega}(v):=\sum_{e\ni v} \omega(e).$$
The irregularity strength of $G$ is in this setting equal to the minimum $k$ such that there exists a weighting $\omega:E\to\{1,2,\ldots,k\}$ with $\sigma_{\omega}(u)\neq \sigma_{\omega}(v)$ for every pair of distinct vertices $u,v$ of $G$. 

A sharp upper bound $s(G)\leq n-1$ for all graphs with $n$ vertices for which the parameter is well defined except $K_3$ was settled in two papers: \cite{Aigner} and~\cite{Nierhoff}, devoted to connected and non-connected cases, respectively, with a $K_{1,n-1}$ star exemplifying its tightness. Nevertheless, much better bound was expected to hold in the case of graphs with larger minimum degree -- a significantly smaller value of the irregularity strength was in particular anticipated for regular graphs. A straightforward counting argument implies a lower bound: $s(G)\geq \lceil\frac{n+d-1}{d}\rceil$ for any $d$-regular graph $G$ with $d\geq 2$. It was conjectured already in 1987 by Faudree and Lehel that this lower bound is not far from optimal. 
\begin{conjecture}[Faudree and Lehel~\cite{Faudree}]\label{FaudreeAndLehelConjecture}
There exists an absolute constant $C$ such that for every $d$-regular graph $G$ with $n$ vertices and $d\geq 2$,
$$s(G)\leq\frac{n}{d}+C.$$
\end{conjecture}
In fact this conjecture was first posed in the form of question by Jacobson (as mentioned by Lehel in~\cite{Lehel}). 
It is this question that triggered extensive studies of the irregularity strength of graphs within the combinatorial community, and resulted in many papers devoted to this graph invariant, see e.g. 
\cite{Aigner,Amar,Amar_Togni,Bohman_Kravitz,Lazebnik,Dinitz,Ebert,Ebert2,Faudree2,Faudree,Ferrara,Frieze,Gyarfas,Jendrol_Tkac,KalKarPf,MajerskiPrzybylo2,Nierhoff,Przybylo,irreg_str2,Togni} 
(or \cite{Lehel} for a survey), giving also rise to many related concepts
\cite{Louigi30,Louigi2,Louigi,AnhKalPrz,Baca,BarGrNiw,LocalIrreg_1,BensmailMerkerThomassen,BonamyPrzybylo,FlandrinMPSW,Kalkowski12,KalKarPf_123,123KLT,Majerski_Przybylo,Przybylo_asym_optim2,Przybylo_asym_optim,LocalIrreg_2,Przybylo_CN_1,1234Reg123,Przybylo_CN_2,12Conjecture,PrzybyloWozniakChoos,Seamon123survey,ThoWuZha,Vuckovic_3-multisets,WongZhu23Choos,WongZhuChoos}, 
to mention just a few out of the vastness of problems of this type.
After over 30 years from its formulation, Conjecture~\ref{FaudreeAndLehelConjecture} still remains open.
The first progress towards its solving was accomplished by Faudree and Lehel themselves, who already in~\cite{Faudree} managed to push the upper bound of $n-1$ down to  $n/2+9$ in the case of regular graphs. 
In 2002 Frieze, Gould, Karo\'nski and Pfender applied the probabilistic method to make another big step forward and show among others that $s(G)\leq 10n/d+1$ when $d\leq\lfloor (n\ln n)^{1/4}\rfloor$ and $s(G)\leq 48n/d+1$ for $d\leq\lfloor n^{1/2}\rfloor$. They also achieved similar results but with slightly larger constants in the case of general graphs (with $d$ replaced by $\delta$). Later Cuckler and Lazebnik again exploited probabilistic approach to prove in particular that $s(G)\leq 48n/d+6$ for a $d$-regular graph with $d\geq 10^{4/3}n^{2/3}\ln^{1/3}n$ (and $s(G)\leq 48n/\delta+6$ for graphs with minimum degree $\delta\geq 10n^{3/4}\ln^{1/4}n$). No linear bounds in $n/d$ (or $n/\delta$) valid for entire spectrum of possible degrees were known at that point. The first such upper bounds were provided by Przyby{\l}o in~\cite{Przybylo} and~\cite{irreg_str2}, where it was proved, resp., that  
$s(G)\leq 16\frac{n}{d}+6$ for $d$-regular graphs and $s(G)\leq 112\frac{n}{\delta}+28$ in the genaral case.
Currently the best result of this type is due to Kalkowski, Karo\'nski and Pfender, who made use of a relatively simple deterministic algorithm, refining and adapting its previous versions designed to tackle related subjects (see e.g.~\cite{Kalkowski12} and~\cite{AnhKalPrz}), and proved in~\cite{KalKarPf}  the general upper bound: $s(G)\leq 6\lceil n/\delta\rceil$ for all graphs with $n$ vertices and minimum degree $\delta\geq 1$
which do not contain isolated edges. Their approach was further developed with enhancement of the probabilistic method by Majerski and Przyby{\l}o, who showed in~\cite{MajerskiPrzybylo2} that $s(G)\leq (4+o(1)) n/\delta+4$ for graphs with $\delta\geq n^{0.5}\ln n$. No better results have been achieved with these techniques for regular graphs.

In this paper we prove that the conjecture of Faudree and Lehel holds asymptotically in the cases when $d$ is neither very small nor very close to n. We in particular show that for large enough $n$ and $d\in [\ln^8n,\frac{n}{\ln^3 n}]$, $s(G)\leq \frac{n}{d}(1+\frac{8}{\ln n})$, thereby proving that $s(G) = \frac{n}{d}(1+o(1))$ then.
We  moreover prove the latter to hold already when $d\in [\ln^{1+\varepsilon}n,\frac{n}{\ln^\varepsilon n}]$ where $\varepsilon$ is an arbitrary positive constant, see Corollaries~\ref{CorollaryLnNice} and~\ref{CorollaryLnBest} in the last section. In fact the both results are implied by the following theorem with a slightly less self-evident form.
\begin{theorem}\label{MainTheoremIrregStrAsympt}
For any positive real numbers $b,\varepsilon$, there exists a constant $N$ such that for every $d$-regular graph $G$
with order $n\geq N$ and $d\in\left[\ln^{1+6b+12\varepsilon}n,\frac{n}{\ln^{2b+5\varepsilon}n}\right]$,
$$s(G)\leq \frac{n}{d}\left(1+\frac{8}{\ln^b n}\right).$$
\end{theorem}

\section{Idea of Proof}

Given a $d$-regular graph $G=(V,E)$ with $n=|V|$ large enough and $d\in[\ln^{1+6b+12\varepsilon}n,$ $\frac{n}{\ln^{2b+5\varepsilon}n}]$, we first randomly choose some special sufficiently large subset $U\subseteq V$, yet still containing a small fraction of all vertices.

We then randomly and independently assign a real number $x_v\in[0,1]$ to every vertex $v$ in $V_0:=V\setminus U$ and associate weight $\lceil\frac{n}{d}\rceil$ to each edge $uv\in E(G[V_0])$ with $x_u+x_v\geq 1$ and weight $0$ otherwise. This more or less yields a desired distribution of weighted degrees in $V_0$, but only roughly.

In order to obtain actual precise distinction of the weighted degrees in $V_0$ we carefully choose weights for edges between $U$ and $V_0$. All of these weights are at the same time chosen relatively large (all appropriately larger than $\lceil\frac{n}{d}\rceil$), so that we are certain that the weighted degrees in $U$ are already larger than all those in $V_0$ (this shall be feasible in particular due to the fact that the randomly chosen $U$ shall be relatively small
and most of the edges incident with any vertex in $U$ shall be joining it with $V_0$).

It is then sufficient to distinguish weighted degrees in $U$ via appropriate choice of weights for the edges of $G[U]$. 
For this goal we shall use an adaptation of the algorithm of Kalkowski, Karo\'nski and Pfender from~\cite{KalKarPf}   (applicable to prove the general $6\lceil n/\delta\rceil$ upper bound).
However, as in our randomly chosen $G[U]$ the proportion of the number of vertices to its minimum degree shall still be close to $\frac{n}{d}$ (as in $G$ itself), we shall have to be extra careful priory while choosing the weights for the edges between $U$ and $V_0$. Namely, we shall guarantee within that process that $U$ shall partition into $7$ subsets $U_1,U_2,\ldots,U_7$ with increasing weighted degrees and no possible conflicts between vertices from distinct subsets. This shall admit an adaptation of the mentioned algorithm with weights reduced up to at most  $\frac{n}{d}+1$.

Within analysis of our random process we shall in particular use the Chernoff Bound. The following its version can be found e.g. in~\cite{JansonLuczakRucinski} (Th. 2.1, page 26).
\begin{theorem}[\textbf{Chernoff Bound}]\label{ChernofBoundTh}
For any $0\leq t\leq np$,
$$\mathbf{Pr}({\rm BIN}(n,p)>np+t)<e^{-\frac{t^2}{3np}}~~~~{and}~~~~\mathbf{Pr}({\rm BIN}(n,p)<np-t)<e^{-\frac{t^2}{2np}}$$ 
where ${\rm BIN}(n,p)$ is the sum of $n$ independent Bernoulli variables, each equal to $1$ with probability $p$ and $0$ otherwise.
\end{theorem}

\section{Proof of Theorem \ref{MainTheoremIrregStrAsympt}}

Let $b,\varepsilon$ be two arbitrarily chosen and fixed positive real numbers.

Let $G=(V,E)$ be a $d$-regular graph with $n=|V|$ and $d\in\left[\ln^{1+6b+12\varepsilon}n,\frac{n}{\ln^{2b+5\varepsilon}n}\right]$. We shall not specify $N$, assuming whenever needed that $n$ is sufficiently large so that explicit inequalities below hold. We shall prove that then 
$$s(G)\leq \frac{n}{d}\left(1+\frac{8}{\ln^b n}\right).$$

\subsection{Random Vertex Partition}

We first observe that we may fix a
partition of $V$ into a smaller part $U$ and a larger part $V_0=V\setminus U$ and a subpartition of $U$
into roughly equal seven parts with proportional distributions of neighbours of every vertex of $G$ between these sets.

\begin{observation}\label{UV-decomposition}
There is a subset $U$ of $V$ and its partition $U=U_1\cup U_2\cup \ldots \cup U_7$ such that for every $v\in V$ and $i\in\{1,2,\ldots,7\}$,
\begin{itemize}
\item[$(1^\circ)$] $~~\left||U_i|-\frac{1}{7}\frac{n}{\ln^{b+\varepsilon} n}\right| \leq \frac{1}{7} \frac{n}{\ln^{2b+4\varepsilon}n}$, 
~~~~~~hence $~~~\left||U|-\frac{n}{\ln^{b+\varepsilon} n}\right| \leq  \frac{n}{\ln^{2b+4\varepsilon}n}$;
\item[$(2^\circ)$] $~~\left|d_{U_i}(v) - \frac{1}{7}\frac{d}{\ln^{b+\varepsilon} n}\right| \leq \frac{1}{7}\frac{d}{\ln^{2b+4\varepsilon} n}$, 
~~~hence $~~~\left|d_{U}(v) - \frac{d}{\ln^{b+\varepsilon} n}\right| \leq \frac{d}{\ln^{2b+4\varepsilon} n}$.
\end{itemize}
\end{observation}

\begin{pf}
We place every vertex $v$ of $G$ in $U$ independently with probability $\frac{1}{\ln^{b+\varepsilon} n}$
and independently we randomly and equiprobably assign an integer $i_v\in\{1,2,\ldots,7\}$ to every vertex $v\in V$.
For every $i=1,2,\ldots,7$, we define the set $U_i$ as the set of vertices $v\in U$ with $i_v=i$.
Then $E(|U_i|)=\frac{1}{7}\frac{n}{\ln^{b+\varepsilon} n}$ and $E(d_{U_i}(v)) = \frac{1}{7}\frac{d}{\ln^{b+\varepsilon} n}$, and by the Chernoff Bound:
$$\mathbf{Pr}\left(\left||U_i|-\frac{1}{7}\frac{n}{\ln^{b+\varepsilon} n}\right| > \frac{1}{7} \frac{n}{\ln^{2b+4\varepsilon}n}\right) < 2e^{-\frac{n}{3\cdot 7\ln^{3b+7\varepsilon}n}} < \frac{1}{14},$$
$$\textbf{Pr}\left( \left|d_{U_i}(v) - \frac{1}{7}\frac{d}{\ln^{b+\varepsilon} n}\right| > \frac{1}{7}\frac{d}{\ln^{2b+4\varepsilon} n}\right) 
< 2e^{-\frac{d}{3\cdot 7\ln^{3b+7\varepsilon}n}} 
\leq 2\left(e^{\ln n}\right)^{-\frac{\ln^{3b+5\varepsilon}n}{21}} 
< \frac{1}{14n}.$$
Therefore, the probability that $(1^\circ)$ does not hold for some $i\in\{1,2,\ldots,7\}$ or $(2^\circ)$ does not hold for some $v\in V$ and $i\in\{1,2,\ldots,7\}$ is less than 
$$7\cdot \frac{1}{14} + n\cdot 7\cdot \frac{1}{14n} = 1,$$
hence there is a choice of $U,U_1,U_2,\ldots,U_7$ fulfilling all our requirements.
\qed
\end{pf}

Denote 
$$V_0:=V\setminus U, ~~~~~~~~n_0:=|V_0|, ~~~~~~~~G_0:=G[V_0]$$ 
and 
$$d_0(v):=d_{V_0}(v)$$
 for every $v\in V$.
Then by~$(1^\circ)$, 
\begin{equation}\label{n_0Bounds}
n\left(1 - \frac{1}{\ln^{b+\varepsilon} n} - \frac{1}{\ln^{2b+4\varepsilon}n}\right) 
\leq n_0 \leq
n\left(1 - \frac{1}{\ln^{b+\varepsilon} n} + \frac{1}{\ln^{2b+4\varepsilon}n}\right)
\end{equation}
and by~$(2^\circ)$ for each $v\in V$:
\begin{equation}\label{d_0Bounds}
d\left(1 - \frac{1}{\ln^{b+\varepsilon} n} - \frac{1}{\ln^{2b+4\varepsilon}n}\right) 
\leq d_0(v) \leq
d\left(1 - \frac{1}{\ln^{b+\varepsilon} n} + \frac{1}{\ln^{2b+4\varepsilon}n}\right).
\end{equation}

\subsection{Random Labeling of the Vertices in $V_0$}

We now aim at assigning weights to all edges with at least one end in $V_0$ so that the obtained weighted degrees in  $V_0$ form an arithmetic progression with step size $1$ (i.e., they are consecutive integers). For this goal we first randomly and independently  choose for every vertex $v$ of $V_0$ a real number $x_v\in[0,1]$ with uniform probability distribution (that is a realization of the random variable $X_v\sim U(0,1)$ associated to $v$). These $x_v$'s shall roughly refer to the positions of vertices in the arithmetic progression and define a natural ordering of the vertices (provided that all $x_v$'s are distinct).

Denote 
 $$L_v:=\{u\in V_0: x_u < x_v\},~~~~~R_v:=\{u\in N_{G_0}(v): x_u\geq 1-x_v\}.$$
In fact $|L_v|$ shall correspond to the number of vertices preceding $v$ in the ordering (and thus to its position within it), while $R_v$ shall define the ends of the edges incident with $v$ in $G_0$ which shall receive weight $\lceil\frac{n}{d}\rceil$ (the remaining ones shall temporarily be weighted $0$).

\begin{observation}\label{GoodX_vDistribution}
With positive probability all $x_v$'s are  pairwise distinct and for every vertex $v\in V_0$ with $x_v=x$:  
\begin{itemize}
\item[$(3^\circ)$]  ~~if~~~ $x\geq \frac{1}{\ln^{2b+3\varepsilon}n}$, ~~~then~~~ $\left||L_v|-x(n_0-1)\right|\leq \frac{x(n_0-1)}{\ln^{2b+4\varepsilon} n}$;
\item[$(4^\circ)$]  ~~if~~~ $x < \frac{1}{\ln^{2b+3\varepsilon}n}$, ~~~then~~~ $|L_v| \leq \frac{n_0-1}{\ln^{2b+3\varepsilon}n} + \frac{n_0-1}{\ln^{4b+7\varepsilon}n}$;
\item[$(5^\circ)$]  ~~if~~~ $x\geq \frac{1}{\ln^{2b+3\varepsilon}n}$, ~~~then~~~ $\left||R_v|-xd_{0}(v)\right|\leq \frac{xd_{0}(v)}{\ln^{2b+4\varepsilon} n}$;
\item[$(6^\circ)$]  ~~if~~~ $x < \frac{1}{\ln^{2b+3\varepsilon}n}$, ~~~then~~~ $|R_v|\leq 
\frac{d_{0}(v)}{\ln^{2b+3\varepsilon} n} + \frac{d_{0}(v)}{\ln^{4b+7\varepsilon} n}$.
\end{itemize}
\end{observation}

\begin{pf}
For every $v\in V_0$ with $x_v=x$ we denote the following events: 
\begin{eqnarray}
A_v: && \left(x\geq \frac{1}{\ln^{2b+3\varepsilon}n}~~\wedge~~\left||L_v|-x(n_0-1)\right| > \frac{x(n_0-1)}{\ln^{2b+4\varepsilon} n}\right) \nonumber\\
&&~~\vee~~
\left(x < \frac{1}{\ln^{2b+3\varepsilon}n}~\wedge~|L_v| > \frac{n_0-1}{\ln^{2b+3\varepsilon}n} + \frac{n_0-1}{\ln^{4b+7\varepsilon}n}\right);\nonumber
\end{eqnarray}
\begin{eqnarray}
B_v: && \left(x\geq \frac{1}{\ln^{2b+3\varepsilon}n}~\wedge~\left||R_v|-xd_{0}(v)\right| > \frac{xd_{0}(v)}{\ln^{2b+4\varepsilon} n}\right) \nonumber\\
&&~~\vee~~\left(x < \frac{1}{\ln^{2b+3\varepsilon}n}~\wedge~|R_v| > 
\frac{d_{0}(v)}{\ln^{2b+3\varepsilon} n} + \frac{d_{0}(v)}{\ln^{4b+7\varepsilon} n}\right).\nonumber
\end{eqnarray}

We note that for any $x\geq \frac{1}{\ln^{2b+3\varepsilon}n}$ 
(as $v\notin L_v$, while each remaining $u\in V_0$ lands in $L_v$ independently with probability $x$),
by the Chernoff Bound and~(\ref{n_0Bounds}):
\begin{eqnarray}
\mathbf{Pr} \left(A_v|X_v=x\right) &=& 
\mathbf{Pr} \left(\left|{\rm BIN}\left(n_0-1,x\right)-x(n_0-1)\right| > \frac{x(n_0-1)}{\ln^{2b+4\varepsilon} n}\right) \nonumber\\
&<& 2e^{-\frac{x(n_0-1)}{3\ln^{4b+8\varepsilon}n}} \leq 2e^{-\frac{n_0-1}{3\ln^{6b+11\varepsilon}n}} < 2e^{-\frac{n}{4\ln^{6b+11\varepsilon}n}}  < \frac{1}{2n}. \label{PrAv1}
\end{eqnarray}
By inequalities in~(\ref{PrAv1}) above (with $x = \frac{1}{\ln^{2b+3\varepsilon}n}$), for any 
$x < \frac{1}{\ln^{2b+3\varepsilon}n}$ we similarly have:
\begin{eqnarray}
\mathbf{Pr} \left(A_v|X_v=x\right) 
&=& \mathbf{Pr} \left({\rm BIN}\left(n_0-1,x\right) > \frac{n_0-1}{\ln^{2b+3\varepsilon}n} + \frac{n_0-1}{\ln^{4b+7\varepsilon}n}\right) \nonumber\\
&\leq& \mathbf{Pr} \left(\left|{\rm BIN}\left(n_0-1,\frac{1}{\ln^{2b+3\varepsilon}n}\right)-\frac{n_0-1}{\ln^{2b+3\varepsilon}n}\right| > \frac{\frac{n_0-1}{\ln^{2b+3\varepsilon}n}}{\ln^{2b+4\varepsilon} n}\right) < \frac{1}{2n}. \label{PrAv2}
\end{eqnarray}

Analogously, for any $x\geq \frac{1}{\ln^{2b+3\varepsilon}n}$, by~(\ref{d_0Bounds}) and the assumption that 
$d\geq\ln^{1+6b+12\varepsilon}n$:
\begin{eqnarray}
\mathbf{Pr} \left(B_v|X_v=x\right) &=& 
\mathbf{Pr} \left(\left|{\rm BIN}\left(d_0(v),x\right)-xd_0(v)\right| > \frac{xd_0(v)}{\ln^{2b+4\varepsilon} n}\right) \nonumber\\
&<& 2e^{-\frac{xd_0(v)}{3\ln^{4b+8\varepsilon}n}} \leq 2e^{-\frac{d_0(v)}{3\ln^{6b+11\varepsilon}n}} < 2e^{-\frac{d}{4\ln^{6b+11\varepsilon}n}}  < \frac{1}{2n}. \label{PrBv1}
\end{eqnarray}
By inequalities in~(\ref{PrBv1}) above, for any $x < \frac{1}{\ln^{2b+3\varepsilon}n}$ we in turn have:
\begin{equation}\label{PrBv2}
\mathbf{Pr} \left(B_v|X_v=x\right) \leq \mathbf{Pr} \left(\left|{\rm BIN}\left(d_0(v),\frac{1}{\ln^{2b+3\varepsilon}n}\right)-
\frac{1}{\ln^{2b+3\varepsilon}n}d_0(v)\right| > \frac{\frac{1}{\ln^{2b+3\varepsilon}n}d_0(v)}{\ln^{2b+4\varepsilon} n}\right) < \frac{1}{2n}.
\end{equation}

By~(\ref{PrAv1}) and~(\ref{PrAv2}), 
$$\mathbf{Pr}(A_v) < \int_0^1\frac{1}{2n} dx = \frac{1}{2n},$$
and by~(\ref{PrBv1}) and~(\ref{PrBv2}), 
$$\mathbf{Pr}(B_v) < \int_0^1\frac{1}{2n} dx = \frac{1}{2n}.$$
Therefore,
$$\mathbf{Pr} \left(\bigcap_{v\in V_0}\left(\overline{A_v}\cap \overline{B_v}\right)\right) = 1 - \mathbf{Pr}\left(\bigcup_{v\in V_0} A_v\cup B_v\right)> 1 - 2n\cdot \frac{1}{2n} = 0,$$
and the thesis follows (as the probability that all $x_v$ are pairwise distinct is obviously $1$).
\qed
\end{pf}

\subsection{Weighted Degrees in $V_0$}

We fix any choice of $x_v$, $v\in V_0$ consistent with Observation~\ref{GoodX_vDistribution} above.

Now to every edge $uv\in E(G_0)$ we assign an initial weight:
$$\omega_0(uv)=\left\{\begin{array}{rcl}
0~~, & {\rm if} & x_u+x_v< 1;\\
\left\lceil\frac{n}{d}\right\rceil,  & {\rm if} & x_u+x_v\geq 1.
\end{array}\right.$$
For every $e\in E[V_0,U_i]$ (the set of edges between $V_0$ and $U_i$) with $i\in\{1,2,\ldots,7\}$ we set:
\begin{equation}\label{Omega0V0U}
\omega_0(e)=\left\lceil\frac{n}{d}\right\rceil + i \left\lceil\frac{n}{d\ln^b n}\right\rceil,
\end{equation}
and for every $e\subset U$ we temporarily set 
$$\omega_0(e)=0.$$

Therefore, $\omega_0(uv) = \left\lceil\frac{n}{d}\right\rceil$ if and only if $x_v \in R_u$ (or equivalently $x_u \in R_v$).

Note that for every $v\in V_0$ we thus have: 
\begin{equation}\label{Omega_0FirstEstimation}
\sigma_{\omega_0} (v) = |R_v|\left\lceil\frac{n}{d}\right\rceil + \sum_{i=1}^7 d_{U_i}(v)\left(\left\lceil\frac{n}{d}\right\rceil + i \left\lceil\frac{n}{d\ln^b n}\right\rceil\right).
\end{equation}

Set
\begin{equation}\label{B_0Definition}
B_0 = \left\lceil\frac{n}{\ln^{b+\varepsilon} n}\right\rceil + 4\left\lceil\frac{n}{\ln^{2b+\varepsilon} n}\right\rceil 
+ 2\left\lceil\frac{n}{\ln^{2b+3\varepsilon} n}\right\rceil,
\end{equation}

\begin{equation}\label{NDefinition}
N=\left\lceil\frac{n}{\ln^{2b+2\varepsilon} n}\right\rceil - 2\left\lceil\frac{n}{\ln^{3b+5\varepsilon} n}\right\rceil.
\end{equation}

By~(\ref{Omega_0FirstEstimation}) and $(2^\circ)$, (\ref{B_0Definition}), (\ref{NDefinition}), we thus have for every given $v\in V_0$ (and $n$ sufficiently large):
\begin{eqnarray}
\sigma_{\omega_0} (v) &\geq& |R_v|\frac{n}{d} + 
\sum_{i=1}^7 \frac{1}{7} \left(\frac{d}{\ln^{b+\varepsilon}n} - \frac{d}{\ln^{2b+4\varepsilon}n}\right)
\left(\frac{n}{d} + i \frac{n}{d\ln^b n}\right) \nonumber\\
&=& |R_v|\frac{n}{d} + 
\frac{n}{\ln^{b+\varepsilon}n} + 4\frac{n}{\ln^{2b+\varepsilon}n} - \frac{n}{\ln^{2b+4\varepsilon}n} - 4\frac{n}{\ln^{3b+4\varepsilon}n}
\nonumber\\
&=& |R_v|\frac{n}{d} + 
\left(\frac{n}{\ln^{b+\varepsilon}n} + 1\right) + 4\left(\frac{n}{\ln^{2b+\varepsilon}n} + 1\right) + 2\left(\frac{n}{\ln^{2b+3\varepsilon}n} + 1\right) \nonumber\\ 
&& - \frac{n}{\ln^{2b+2\varepsilon}n} + 2\left(\frac{n}{\ln^{3b+5\varepsilon}n} + 1\right) 
+ \left(\frac{n}{\ln^{2b+3\varepsilon} n} + \frac{n}{\ln^{4b+7\varepsilon}n}\right) \nonumber\\
&&
+ \left( \frac{n}{\ln^{2b+2\varepsilon}n} - 3\frac{n}{\ln^{2b+3\varepsilon}n}
- \frac{n}{\ln^{2b+4\varepsilon}n} - 4\frac{n}{\ln^{3b+4\varepsilon}n} - 2\frac{n}{\ln^{3b+5\varepsilon}n} - \frac{n}{\ln^{4b+7\varepsilon}n}\ - 9\right)
\nonumber\\
&>& |R_v|\frac{n}{d} +  B_0 - N
+ \left(\frac{n}{\ln^{2b+3\varepsilon} n} + \frac{n}{\ln^{4b+7\varepsilon}n}\right). \label{Sigma0Larger1}
\end{eqnarray}
Thus if $x_v=x < \frac{1}{\ln^{2b+3\varepsilon}n}$, by (\ref{Sigma0Larger1}) and $(4^\circ)$:
\begin{eqnarray}
\sigma_{\omega_0} (v) &>& 
B_0 - N
+ \left(\frac{n}{\ln^{2b+3\varepsilon} n} + \frac{n}{\ln^{4b+7\varepsilon}n}\right) \geq B_0 + |L_v| - N. \label{Sigma0Larger2}
\end{eqnarray}
Analogously, if $x_v=x \geq \frac{1}{\ln^{2b+3\varepsilon}n}$, by (\ref{Sigma0Larger1}) and $(5^\circ)$, (\ref{d_0Bounds}), (\ref{n_0Bounds}), $(3^\circ)$:
\begin{eqnarray}
\sigma_{\omega_0} (v) &>& 
 x d_0(v)
 \left(1 - \frac{1}{\ln^{2b+4\varepsilon} n}\right)\frac{n}{d} 
+ B_0 - N
+ \left(\frac{n}{\ln^{2b+3\varepsilon} n} + \frac{n}{\ln^{4b+7\varepsilon}n}\right) \nonumber\\ 
&\geq& 
 x d\left(1 - \frac{1}{\ln^{b+\varepsilon} n} - \frac{1}{\ln^{2b+4\varepsilon}n}\right) 
 \left(1 - \frac{1}{\ln^{2b+4\varepsilon} n}\right)\frac{n}{d} 
+ B_0 - N
+ \left(\frac{n}{\ln^{2b+3\varepsilon} n} + \frac{n}{\ln^{4b+7\varepsilon}n}\right) \nonumber\\ 
&>&  xn - \frac{xn}{\ln^{b+\varepsilon}n} +
\frac{1}{2} \frac{n}{\ln^{2b+3\varepsilon}n}  + B_0 - N  \nonumber\\ 
&>&  x n  \left(1 - \frac{1}{\ln^{b+\varepsilon} n} + \frac{1}{\ln^{2b+4\varepsilon} n}\right)  
\left(1 + \frac{1}{\ln^{2b+4\varepsilon} n}\right) + B_0 - N  \nonumber\\
&\geq&  x n_0  
\left(1 + \frac{1}{\ln^{2b+4\varepsilon} n}\right) + B_0 - N  \nonumber\\
&>&  B_0 + |L_v| - N. \label{Sigma0Larger3}
\end{eqnarray}

On the other hand, by~(\ref{Omega_0FirstEstimation}) and $(2^\circ)$, (\ref{B_0Definition}), for every given $v\in V_0$ (and $n$ sufficiently large):
\begin{eqnarray}
\sigma_{\omega_0} (v) &\leq& |R_v|\left(\frac{n}{d} + 1\right) 
+ \sum_{i=1}^7 \frac{1}{7} \left(\frac{d}{\ln^{b+\varepsilon}n} + \frac{d}{\ln^{2b+4\varepsilon}n}\right)
\left(\frac{n}{d} + i \frac{n}{d\ln^b n} + 8\right) \nonumber\\
&<& |R_v|\frac{n}{d} + d
+ \sum_{i=1}^7 \frac{1}{7}\frac{d}{\ln^{b+\varepsilon}n}\left(\frac{n}{d} + i \frac{n}{d\ln^b n}\right)
+ \sum_{i=1}^7 \frac{1}{7}\frac{d}{\ln^{2b+4\varepsilon}n} \cdot 2\frac{n}{d} + 9 \frac{d}{\ln^{b+\varepsilon} n}  \nonumber\\
&=&   |R_v|\frac{n}{d} 
+ \frac{n}{\ln^{b+\varepsilon}n} + 4\frac{n}{\ln^{2b+\varepsilon}n} + 2\frac{n}{\ln^{2b+4\varepsilon}n} 
+ d + 9 \frac{d}{\ln^{b+\varepsilon} n}\nonumber\\
&<&   |R_v|\frac{n}{d} 
+ \frac{n}{\ln^{b+\varepsilon}n} + 4\frac{n}{\ln^{2b+\varepsilon}n} + 3\frac{n}{\ln^{2b+4\varepsilon}n} \nonumber\\
&<&   |R_v|\frac{n}{d} + B_0
- \frac{3}{2} \frac{n}{\ln^{2b+3\varepsilon}n}. \label{Sigma0Smaller1}
\end{eqnarray}
Thus if 
$x_v=x < \frac{1}{\ln^{2b+3\varepsilon}n}$, by (\ref{Sigma0Smaller1}) and $(6^\circ)$:
\begin{eqnarray}
\sigma_{\omega_0} (v) &<& 
\left( \frac{d}{\ln^{2b+3\varepsilon} n} + \frac{d}{\ln^{4b+7\varepsilon} n}\right)\frac{n}{d}
+B_0 
- \frac{3}{2} \frac{n}{\ln^{2b+3\varepsilon}n} <
B_0  \leq B_0 + |L_v|. \label{Sigma0Smaller2}
\end{eqnarray}
Analogously, if $x_v=x \geq \frac{1}{\ln^{2b+3\varepsilon}n}$, by (\ref{Sigma0Smaller1}) and $(5^\circ)$, (\ref{d_0Bounds}), (\ref{n_0Bounds}), $(3^\circ)$:
\begin{eqnarray}
\sigma_{\omega_0} (v) &<& 
 x d_0(v)
 \left(1 + \frac{1}{\ln^{2b+4\varepsilon} n}\right)\frac{n}{d} 
+ B_0
- \frac{3}{2} \frac{n}{\ln^{2b+3\varepsilon}n}  \nonumber\\
&\leq&  xd \left(1 - \frac{1}{\ln^{b+\varepsilon} n} + \frac{1}{\ln^{2b+4\varepsilon}n}\right)  \left(1 + \frac{1}{\ln^{2b+4\varepsilon} n}\right)\frac{n}{d}
+ B_0
- \frac{3}{2} \frac{n}{\ln^{2b+3\varepsilon}n} \nonumber\\
&<& xn - \frac{xn}{\ln^{b+\varepsilon} n} + B_0 - 
\frac{n}{\ln^{2b+3\varepsilon}n} \nonumber\\
&<&  x n  \left(1 - \frac{1}{\ln^{b+\varepsilon} n} - \frac{1}{\ln^{2b+4\varepsilon} n} - \frac{1}{n}\right)  
\left(1 - \frac{1}{\ln^{2b+4\varepsilon} n}\right) + B_0 \nonumber\\
&\leq&  x (n_0-1) 
\left(1 - \frac{1}{\ln^{2b+4\varepsilon} n}\right) + B_0 \nonumber\\
&\leq& B_0 + |L_v|. \label{Sigma0Smaller3}
\end{eqnarray}

Fix the (unique) ordering $v_1,v_2,\ldots,v_{n_0}$ of the vertices in $V_0$ so that 
$$x_{v_i}\leq x_{v_j}~~ {\rm whenever}~~ i\leq j.$$ 

Now we admit some changes on the edges $e\in E[V_0,U]$ -- namely for each of these we admit adding any integer
\begin{equation}\label{OmegaPrimInterval}
\omega'(e)\in \left[0,\left\lceil\frac{n}{d\ln^{b+\varepsilon} n}\right\rceil\right]
\end{equation}
to its current weight. For definiteness, we set $\omega'(e)=0$ for the remaining edges of $G$.

As a result we may add any integer in 
$$\left[0,d_U(v)\left\lceil\frac{n}{d\ln^{b+\varepsilon} n}\right\rceil\right]$$
to $\sigma_{\omega_0}(v)$ of any $v\in V_0$. Thus if we set 
$$\omega_1:=\omega_0+\omega',$$
as $(2^\circ)$ implies that
$$d_U(v)\left\lceil\frac{n}{d\ln^{b+\varepsilon} n}\right\rceil 
\geq d\left(\frac{1}{\ln^{b+\varepsilon} n}-\frac{1}{\ln^{2b+4\varepsilon} n}\right) \frac{n}{d\ln^{b+\varepsilon} n} > N,$$
by~(\ref{Sigma0Larger2}),~(\ref{Sigma0Larger3}), ~(\ref{Sigma0Smaller2}) and~(\ref{Sigma0Smaller3}),   
we can make our choices of $\omega'(e)$ so that~(\ref{OmegaPrimInterval}) holds and 
$$\sigma_{\omega_1}(v_j)=B_0+j$$
for every $j\in\{1,2,\ldots,n_0\}$.

\subsection{Distinguishing Vertices in $U$}

We shall not change weighted degrees of vertices in $V_0$ further on, but we shall increase the weights of (some of) the edges in $U$ in order to adjust the weighted degrees of vertices in $U$. We first observe that their weighted degrees are already larger than those in $V_0$. Note for this aim that for every $v\in V_0$ and $u\in U$, by (\ref{n_0Bounds}), (\ref{B_0Definition}), (\ref{d_0Bounds}) and~(\ref{Omega0V0U}):
\begin{equation}\label{sigma1vu}
\sigma_{\omega_1}(v)\leq B_0+n_0 <  
n + 5\frac{n}{\ln^{2b+\varepsilon} n}
<
d_0(u) \left(\left\lceil\frac{n}{d}\right\rceil +  \left\lceil\frac{n}{d\ln^b n}\right\rceil\right)
\leq \sigma_{\omega_1}(u).
\end{equation}

In order to distinguish vertices in $U$ we shall now use an adaptation of the algorithm of Kalkowski, Karo\'nski and Pfender from~\cite{KalKarPf}. Within this we shall be admitting adding integers:
\begin{equation}\label{Omega''Limits}
\omega''(e)\in\left[0,\frac{n}{d}\right] 
\end{equation}
to the weight of every edge $e\subset U$. The almost final weighting of $G$ shall be defined as 
\begin{equation}\label{Omega2asSumOmega1Omega''}
\omega_2=\omega_1+\omega''
\end{equation}
(where we set $w''(e)=0$ for every edge $e\in E\setminus E(G[U])$).
Note in particular that by~(\ref{sigma1vu}), for every $v\in V_0$ and $u\in U$ we shall (still) have:
\begin{equation}\label{DistinctionBetweenUV0}
\sigma_{\omega_2}(v) = \sigma_{\omega_1}(v) <  \sigma_{\omega_1}(u) \leq \sigma_{\omega_2}(u).
\end{equation}

Moreover, so that the algorithm could work effectively, we observe now that we shall be able to focus within its execution on distinguishing merely the vertices within the same $U_i$'s, as by our construction, for every $u\in U_i$ and $v\in U_{i+1}$ with $i\in\{1,2,3,4,5,6\}$, by~(\ref{Omega0V0U}), (\ref{OmegaPrimInterval}), (\ref{Omega''Limits}) and~(\ref{d_0Bounds}) we shall have:
\begin{eqnarray}
\sigma_{\omega_2}(u) 
&\leq& d_0(u) \left(\left\lceil\frac{n}{d}\right\rceil + i \left\lceil\frac{n}{d\ln^b n}\right\rceil + \left\lceil\frac{n}{d\ln^{b+\varepsilon} n}\right\rceil\right)  + d_U(u)\left\lceil\frac{n}{d}\right\rceil  \nonumber\\
&=& d\left\lceil\frac{n}{d}\right\rceil  + d_0(u) \left(i \left\lceil\frac{n}{d\ln^b n}\right\rceil + \left\lceil\frac{n}{d\ln^{b+\varepsilon} n}\right\rceil\right)   \nonumber\\
&<& n +  i\frac{n}{\ln^b n} +  \frac{n}{\ln^{b+\varepsilon} n} \nonumber\\
&<& n +  i\frac{n}{\ln^b n} +  \frac{n}{\ln^b n} - 2\frac{n}{\ln^{b+\varepsilon} n} \nonumber\\
&<& d_0(v) \left(\left\lceil\frac{n}{d}\right\rceil + (i+1) \left\lceil\frac{n}{d\ln^b n}\right\rceil\right) \nonumber\\
&\leq& \sigma_{\omega_2}(v). \label{DistinctionBetweenUis}
\end{eqnarray}

Initially we set 
$$\omega''(e)=\left\lfloor\frac{n}{3d}\right\rfloor$$
for every edge $e$ of $G[U]$. This shall be modified by the algorithm (while $\omega_2$ shall always refer to an up-to-date value of $\omega''$ below; cf.~(\ref{Omega2asSumOmega1Omega''})).

Within the algorithm we then analyze one component of $G[U]$ after another in arbitrary order, and in each consecutive component of $G[U]$, say $G'=(V',E')$, we arrange its vertices into a sequence $u_1,u_2,\ldots,u_{n'}$ so that each of these vertices, except the last one has a \emph{forward edge}, i.e. an edge joining it in $G'$ with a vertex later in the ordering, which we shall call a \emph{forward neighbour} of this vertex (we may use e.g. a reversed BFS ordering for this goal). We define \emph{backward edges} and \emph{neighbours} of a vertex correspondingly.

We analyze all vertices in the sequence one by one consistently with the fixed ordering in $G'$, and for any currently analyzed vertex $v$, we admit adding any integer in 
\begin{equation}\label{ForwardChanges}
\left[0,\left\lfloor\frac{n}{3d}\right\rfloor\right]
\end{equation}
to the weight of its every forward edge and adding one of the three numbers from the set 
\begin{equation}\label{ThreeModifOptions}
\left\{-\left\lfloor\frac{n}{3d}\right\rfloor,0,\left\lfloor\frac{n}{3d}\right\rfloor\right\}
\end{equation}
to the weights of its backward edges (thus every edge of $G[U]$ shall be modified at most twice within the algorithm). However, not all the mentioned three options in~(\ref{ThreeModifOptions}) shall be available for a given backward edge. Namely, the moment a given vertex $v$ is analyzed, we assign to it a suitable two-element set $\Sigma_v$ belonging to the family:
$$S=\left\{\left\{2\lambda \left\lfloor\frac{n}{3d}\right\rfloor +a,(2\lambda+1)\left\lfloor\frac{n}{3d}\right\rfloor+a\right\}: \lambda\in\mathbb{Z}, a\in\left\{0,1,2,\ldots,\left\lfloor\frac{n}{3d}\right\rfloor-1\right\}\right\}$$
(note $S$ \emph{partitions} $\mathbb{Z}$ into two-element sets)
and perform admitted modifications of the weights of its forward and backward edges so that the obtained weighted degree $\sigma_{\omega_2}(v)$ of $v$ belongs to $\Sigma_v$. Once $\Sigma_v$ is fixed, $\sigma_{\omega_2}(v)$ is required to be its element (i.e. it may take only two distinct values from this point till the end of the construction). Therefore, if $uv$ is a backward edge of $v$ and $u=\min\Sigma_u$, then we may increase the weight of $uv$ by $\left\lfloor\frac{n}{3d}\right\rfloor$ (or do not alter it at all), while otherwise we may decrease it by $\left\lfloor\frac{n}{3d}\right\rfloor$. Such admitted changes of the weights of the backward edges of $v$ along with the admitted alterations concerning forward edges (cf.~(\ref{ForwardChanges})) guarantee by~$(2^\circ)$ at least
\begin{equation}
\label{NoOptionSima2}
d_U(v)\cdot \left\lfloor\frac{n}{3d}\right\rfloor+1 > 2 \cdot\frac{n}{7\ln^{b+\varepsilon}n}\left(1+\frac{1}{\ln^{b+3\varepsilon}n}\right)
\end{equation}
options for  $\sigma_{\omega_2}(v)$ for every currently analyzed vertex $v$ except for the last vertex of the component 
(which is the only vertex of $G'$ without a forward edge). 
Suppose $v\in U_i$ for some $i\in\{1,2,\ldots,7\}$ and $v\notin\{u_{n'-1},u_{n'}\}$.
Then we may perform the admissible changes so that $\sigma_{\omega_2}(v)$ does not belong to $\Sigma_u$ for every vertex $u\in U_i$ (which has $\Sigma_u$ already fixed), as by~(\ref{NoOptionSima2}) and~$(1^\circ)$, $2|U_i|<  
d_U(v)\cdot \left\lfloor\frac{n}{3d}\right\rfloor+1$. Then we fix as $\Sigma_v$ the only set in $S$ which contains the obtained $\sigma_{\omega_2}(v)$.

Now suppose we have already analyzed $u_1,u_2,\ldots,u_{n'-2}$, hence we are left with just two vertices to by analyzed in the given component $G'$. Note that by our construction and~(\ref{DistinctionBetweenUis}), each $\Sigma_{u_j}$ with $j=1,2,\ldots,n'-2$ is disjoint with $\Sigma_u$ for every vertex $u\neq u_j$ which has already been analyzed. Let $L$ be the set consisting of all last and last but one vertices in the components of $G[U]$, hence $L\cap V'=\{u_{n'-1},u_{n'}\}$. The vertices in $L$ shall be admitted to have weighted degrees belonging to some previously 
fixed $\Sigma_u$ (yet different from $\sigma_{\omega_2}(u)$). Note that as this consent concerns only vertices in $L$ (and the corresponding sets $\Sigma_u$), the number of sets $\Sigma\in S$ assigned to more than one vertex equals at most twice the number of components of $G[U]$, i.e., by~$(1^\circ)$ and~$(2^\circ)$, less than 
$$2\cdot \frac{|U|}{\delta(G[U])}  < 3\frac{n}{d}$$ 
sets. Denote temporarily the family of all such sets by $S_t$. Therefore, we may change (choose) the current value $\omega''(u_{n'-1}u_{n'}) = \left\lfloor\frac{n}{3d}\right\rfloor$ to some quantity in $\left\{0,1,\ldots,\left\lfloor\frac{n}{3d}\right\rfloor\right\}$ so that at most
$$ \left\lfloor \frac{2\cdot\frac{3n}{d}}{\left\lfloor\frac{n}{3d}\right\rfloor} \right\rfloor < 20 $$
sets in $S_t$ contain elements congruent to the resulting $\sigma_{\omega_2}(u_{n'-1})$ and at the same time at most the same number (i.e. 20) of these sets contain elements congruent to the resulting 
$\sigma_{\omega_2}(u_{n'})$ modulo $\left\lfloor\frac{n}{3d}\right\rfloor$. 
Then the admitted alterations on the backward edges of $u_{n'-1}$ provide $d_U(u_{n'-1})\geq 45$ options (cf.~$(2^\circ)$) for the weighted degree of  $u_{n'-1}$ which form an arithmetic progression with step size $\left\lfloor\frac{n}{3d}\right\rfloor$. Thus at least one of these, say $\sigma'\in \Sigma'\in S$, does not belong to any set in $S_t$ and is neither the largest nor the smallest element of this arithmetic progression and moreover is not the current weighted degree of any already analyzed vertex -- we then perform admissible modifications of the weights of backward edges of $u_{n'-1}$ to obtain $\sigma_{\omega_2}(u_{n'-1})=\sigma'$ not using for this aim a potential at most one edge $u'u_{n'-1}$ with $\sigma_{\omega_2}(u')\in \Sigma'$ (note it is always possible, as $\sigma'$ was chosen not to be the larges nor the smallest element of the mentioned arithmetic progression). We then set $\Sigma_{u_{n'-1}}=\Sigma'$ and note that the obtained  $\sigma_{\omega_2}(u_{n'-1})$ is distinct from weighted degrees of all already analyzed vertices of $G'$. 

Finally we analyze $u_{n'}$. The admitted alterations on the backward edges of $u_{n'}$ distinct from $u_{n'}u_{n'-1}$ and from a possible single edge $u_{n'}u'$ with $\sigma_{\omega_2}(u')\in \Sigma'$ provide $d_U(u_{n'})-1\geq 47$ options (cf.~$(2^\circ)$) for the weighted degree of  $u_{n'}$ which form an arithmetic progression with step size $\left\lfloor\frac{n}{3d}\right\rfloor$. Thus at least one of these, say $\sigma''\in \Sigma''\in S$, does not belong to any set in $S_t\cup\{\Sigma'\}$ and is neither the largest nor the smallest element of this arithmetic progression, and additionally is not the current weighted degree of any already analyzed vertex -- we then perform admissible modifications of the weights of backward edges of $u_{n'}$ distinct from $u_{n'}u_{n'-1}$ and the possible edge $u_{n'}u'$ to obtain $\sigma_{\omega_2}(u_{n'})=\sigma''$ not using for this aim (analogously as previously) a potential at most one edge $u''u_{n'}$ with $\sigma_{\omega_2}(u'')\in \Sigma''$. We then set $\Sigma_{u_{n'}}=\Sigma''$ and note that all the already analyzed vertices of $G'$ have distinct current weighted degrees, which shall not be altered in the further part of the construction.

After analyzing all vertices of $G[U]$ we obtain a weighting $\omega_2$ such that $\sigma_{\omega_2}(u)\neq \sigma_{\omega_2}(v)$ for every $uv\in E$, as by our construction we have for every $e\in E(G[U])$: 
$$\omega''(e)\in\left[0,3\left\lfloor\frac{n}{3d}\right\rfloor\right],$$
hence~(\ref{Omega''Limits}) indeed holds, and thus in particular also~(\ref{DistinctionBetweenUV0}) and~(\ref{DistinctionBetweenUis})  are fulfilled. We then finally set:
$$\sigma_3:=1+\sigma_2,$$
and as $G$ is regular, we obtain that for every $uv\in E$:
$$\sigma_{\omega_3}(u)\neq \sigma_{\omega_3}(v)$$ 
and
$$1\leq \sigma_3(uv) \leq \left\lceil\frac{n}{d}\right\rceil + 7 \left\lceil\frac{n}{d\ln^b n}\right\rceil + \left\lceil\frac{n}{d\ln^{b+\varepsilon} n}\right\rceil +1 \leq \frac{n}{d} + 8\frac{n}{d\ln^b n}.$$
\qed

\section{Concluding Remarks}

By substituting e.g. $b=1$, $\varepsilon=\frac{1}{12}$ in Theorem~\ref{MainTheoremIrregStrAsympt} we obtain the following:
\begin{corollary}\label{CorollaryLnNice}
For every $d$-regular graph with $n$ vertices and $d\in\left[\ln^{8}n,\frac{n}{\ln^3n}\right]$:
$$s(G)\leq \frac{n}{d}\left(1+\frac{8}{\ln n}\right)$$
for $n$ sufficiently large.
\end{corollary}
So that our probabilistic argument works effectively, a poly-logarithmic in $n$ lower bound on $d$ is unfortunately unavoidable. We may however still conclude that $s(G)=(1+o(1))n/d$ for a wider domain of $d$ e.g. by fixing any small yet positive $\varepsilon_0$ and substituting $b=\frac{\varepsilon_0}{18}$ (i.e., $b>\frac{\varepsilon_0}{19}$) and $\varepsilon=\frac{\varepsilon_0}{18}$ in Theorem~\ref{MainTheoremIrregStrAsympt} to obtain: 
\begin{corollary}\label{CorollaryLnBest}
For each fixed $\epsilon_0>0$, for every $d$-regular graph with $n$ vertices and $d\in\left[\ln^{1+\varepsilon_0}n,\frac{n}{\ln^{\varepsilon_0}n}\right]$, if $n$ is sufficiently large,
$$s(G)\leq \frac{n}{d}\left(1+\frac{1}{\ln^{\frac{\varepsilon_0}{19}} n}\right).$$
Hence,
$$s(G)=(1+o(1))\frac{n}{d}.$$
\end{corollary}
We also note that on the other hand the second order term in our upper bounds above could also be greatly improved, but at the cost of narrowing down the interval for $d$. In particular, using the same technique as in the proof of Theorem~\ref{MainTheoremIrregStrAsympt}  one may show that $s(G)\leq\frac{n}{d}(1+\frac{1}{n^\gamma})$ for any $d$-regular graph $G$ with sufficiently large order $n$ and $d\in[n^{\alpha},n^{\beta}]$ for appropriately selected constants  $\alpha,\beta,\gamma\in (0,1)$, $\alpha<\beta$. We omit derails here, as the major goal of this paper was to settle an asymptotic confirmation of Conjecture~\ref{FaudreeAndLehelConjecture}, namely that  $s(G)=(1+o(1))\frac{n}{d}$, for 
possibly widest spectrum of degrees $d$.


\begin{thebibliography}{99}

\bibitem{Louigi30}
L. Addario-Berry, K. Dalal, C. McDiarmid, B.A. Reed, A. Thomason,
\emph{Vertex-Colouring Edge-Weightings}, Combinatorica 27(1) (2007) 1--12.

\bibitem{Louigi2}
L. Addario-Berry, R.E.L. Aldred, K. Dalal, B.A. Reed, \emph{Vertex colouring edge partitions}, 
J. Combin. Theory Ser. B, 94(2) (2005) 237--244.

\bibitem{Louigi}
L. Addario-Berry, K. Dalal, B.A. Reed, \emph{Degree Constrained
Subgraphs},
Discrete Appl. Math. 156(7) (2008) 1168--1174.

\bibitem{Aigner}
M. Aigner, E. Triesch, \emph{Irregular assignments of trees and
forests}, SIAM J. Discrete Math. 3(4) (1990) 439--449.

\bibitem{Amar}
D. Amar, Irregularity strength of regular graphs of large degree, Discrete Math. 114 (1993) 9--17.

\bibitem{Amar_Togni}
D. Amar, O. Togni, Irregularity strength of trees, Discrete Math. 190 (1998)  15--38.

\bibitem{AnhKalPrz}
M. Anholcer, M. Kalkowski, J. Przyby{\l}o, A new upper bound for the total vertex irregularity strength of graphs,
Discrete Math. 309 (2009) 6316--6317.

\bibitem{Baca}
M. Ba\v ca, S. Jendro\v l, M. Miller, J. Ryan, On irregular
total labellings, Discrete Math. 307 (2007) 1378--1388.

\bibitem{BarGrNiw}
T. Bartnicki, J. Grytczuk, S. Niwczyk, \emph{Weight Choosability
of Graphs}, J. Graph Theory 60(3) (2009) 242--256.

\bibitem{LocalIrreg_1}
O. Baudon, J. Bensmail, J. Przyby{\l}o, M. Wo\'zniak, \emph{On decomposing regular graphs into locally irregular subgraphs}, European J. Combin. 49 (2015) 90--104.

\bibitem{BensmailMerkerThomassen}
J. Bensmail, M. Merker, C. Thomassen, \emph{Decomposing graphs into a constant number
of locally irregular subgraphs}, European J. Combin. 60 (2017) 124--134.

\bibitem{Bohman_Kravitz}
T. Bohman, D. Kravitz, On the irregularity strength of trees, J. Graph Theory
45 (2004) 241--254.

\bibitem{BonamyPrzybylo}
M. Bonamy, J. Przyby{\l}o, \emph{On the neighbor sum distinguishing index of planar graphs},
J. Graph Theory 85(3) (2017) 669--690.

\bibitem{ChartrandErdosOellermann}
G. Chartrand, P. Erd\H{o}s, O.R. Oellermann, \emph{How to Define an Irregular Graph},
College Math. J. 19(1) (1988) 36--42.

\bibitem{Chartrand}
G. Chartrand, M.S. Jacobson, J. Lehel, O.R. Oellermann, S. Ruiz, F. Saba,
\emph{Irregular networks}, Congr. Numer. 64 (1988) 197--210.

\bibitem{Lazebnik}
B. Cuckler, F. Lazebnik, \emph{Irregularity Strength of Dense
Graphs}, J. Graph Theory 58(4) (2008) 299--313.

\bibitem{Dinitz}
J.H. Dinitz, D.K. Garnick, A. Gy\'arf\'as, On the irregularity strength of the $m \times n$
grid, J. Graph Theory 16 (1992) 355--374.

\bibitem{Ebert}
G. Ebert, J. Hemmeter, F. Lazebnik, A.J. Woldar, On the irregularity strength
of some graphs, Congr. Numer. 71 (1990) 39--52.

\bibitem{Ebert2}
G. Ebert, J. Hemmeter, F. Lazebnik, A.J. Woldar, On the number of irregular
assignments on a graph, Discrete Math. 93 (1991) 131--142.

\bibitem{Faudree2}
R.J. Faudree, M.S. Jacobson, J. Lehel, R. Schelp, \emph{Irregular
networks, regular graphs and integer matrices with distinct row
and column sums}, Discrete Math. 76 (1989) 223--240.

\bibitem{Faudree}
R.J. Faudree, J. Lehel, \emph{Bound on the irregularity strength
of regular graphs}, Colloq Math Soc Ja\'nos Bolyai, 52,
Combinatorics, Eger North Holland, Amsterdam, (1987), 247--256.

\bibitem{Ferrara}
M. Ferrara, R.J. Gould, M. Karo\'nski, F. Pfender, An iterative approach to graph irregularity strength, Discr. Appl. Math. 158 (2010) 1189--1194.

\bibitem{FlandrinMPSW}
E. Flandrin, A. Marczyk, J. Przyby{\l}o, J-F. Sacl\'e, M. Wo{\'z}niak,
\emph{Neighbor sum distinguishing index},
Graphs Combin. 29(5) (2013) 1329--1336.

\bibitem{Frieze}
A. Frieze, R.J. Gould, M. Karo\'nski, F. Pfender, \emph{On Graph
Irregularity Strength}, J. Graph Theory 41(2) (2002)
120--137.

\bibitem{Gyarfas}
A. Gy\'arf\'as, The irregularity strength of $K_{m,m}$ is
$4$ for odd $m$, Discrete Math. 71 (1998) 273--274.

\bibitem{JansonLuczakRucinski}
S. Janson, T. {\L}uczak, A. Ruci\'nski, \emph{Random Graphs}, Wiley, New York, 2000.

\bibitem{Jendrol_Tkac}
S. Jendrol', M. Tka\v c, Z. Tuza, The irregularity strength and cost of the union of
cliques, Selected papers in honour of Paul Erd\H os on the occasion of his 80th birthday
(Keszthely, 1993), Discrete Math. 150 (1996) 179--186.

\bibitem{Kalkowski12}
M. Kalkowski, \emph{A note on 1,2-Conjecture}, in Ph.D. Thesis, 2009.
Available at:
https://repozytorium.amu.edu.pl/bitstream/10593/445/1\\
/Rozprawa\%20Doktorska\%202010\%20Maciej\%20Kalkowski.pdf

\bibitem{KalKarPf}
M. Kalkowski, M. Karo\'nski, F. Pfender, \emph{A new upper bound for the irregularity strength of graphs},
SIAM J. Discrete Math. 25 (2011) 1319--1321.

\bibitem{KalKarPf_123}
M. Kalkowski, M. Karo\'nski, F. Pfender, \emph{Vertex-coloring edge-weightings: Towards the 1-2-3 conjecture}, J. Combin. Theory Ser. B 100 (2010) 347--349.

\bibitem{123KLT}
M. Karo\'nski, T. \L uczak, A. Thomason, \emph{Edge weights and
vertex colours}, J. Combin. Theory Ser. B 91 (2004) 151--157.

\bibitem{Lehel}
J. Lehel, 
\emph{Facts and quests on degree irregular assignments},
Graph Theory, Combinatorics and Applications, Willey, New York, 1991, 765--782.

\bibitem{MajerskiPrzybylo2}
P. Majerski, J. Przyby{\l}o, \emph{On the irregularity strength of dense graphs},
SIAM J. Discrete Math. 28(1) (2014) 197--205.

\bibitem{Majerski_Przybylo}
P. Majerski, J. Przyby{\l}o, Total vertex irregularity strength of dense graphs. {\em J. Graph Theory} 76(1) (2014) 34--41.

\bibitem{Nierhoff}
T. Nierhoff, \emph{A tight bound on the irregularity strength of
graphs}, SIAM J. Discrete Math. 13(3) (2000) 313--323.

\bibitem{Przybylo_asym_optim2}
J. Przyby{\l}o,
\emph{A note on asymptotically optimal neighbour sum distinguishing colourings},
European J. Combin. 77 (2019) 49--56.

\bibitem{Przybylo_asym_optim}
J. Przyby{\l}o,
\emph{Asymptotically optimal neighbour sum distinguishing colourings of graphs},
Random Structures Algorithms 47 (2015) 776--791.

\bibitem{LocalIrreg_2}
J. Przyby{\l}o, \emph{On decomposing graphs of large minimum degree into locally irregular subgraphs}, Electron. J. Combin. 23(2) (2016) $\sharp$P2.31.

\bibitem{Przybylo}
J. Przyby{\l}o, \emph{Irregularity strength of regular graphs},
Electron. J. Combin. 15(1) (2008) $\sharp$R82.

\bibitem{irreg_str2}
J. Przyby{\l}o, \emph{Linear bound on the irregularity strength and the
total vertex irregularity strength of graphs}, SIAM J.
Discrete Math. 23(1) (2009) 511--516.

\bibitem{Przybylo_CN_1}
J. Przyby{\l}o, \emph{Neighbor distinguishing edge colorings via the Combinatorial Nullstellensatz}, SIAM J. Discrete Math. 27(3) (2013) 1313--1322.

\bibitem{1234Reg123}
J. Przyby{\l}o, \emph{The 1--2--3 Conjecture almost holds for regular graphs}, submitted.

\bibitem{Przybylo_CN_2}
J. Przyby{\l}o, T-L. Wong, \emph{Neighbor distinguishing edge colorings via the Combinatorial Nullstellensatz revisited}, J. Graph Theory 80(4) (2015) 299--312.

\bibitem{12Conjecture}
J. Przyby{\l}o, M. Wo\'zniak, \emph{On a 1,2 Conjecture},
Discrete Math. Theor. Comput. Sci. 12(1) (2010) 101--108.

\bibitem{PrzybyloWozniakChoos}
J. Przyby{\l}o, M. Wo\'zniak, \emph{Total weight choosability of graphs},
Electron. J. Combin. 18(1) (2011) $\sharp$P112.

\bibitem{Seamon123survey}
B. Seamone, \emph{The 1-2-3 Conjecture and related problems: a survey}, Technical report,
available online at http://arxiv.org/abs/1211.5122, 2012.

\bibitem{ThoWuZha} C. Thomassen, Y. Wu, C.Q. Zhang, \emph{The $3$-flow conjecture, factors modulo $k$, and the 1--2--3 conjecture}, J. Combin. Theory Ser. B 121 (2016) 308--325.

\bibitem{Togni}
O. Togni, Irregularity strength and compound graphs, Discrete Math. 218 (2000) 235--243.

\bibitem{Vuckovic_3-multisets}
B. Vu\v{c}kovi\'c, \emph{Multi-set neighbor distinguishing 3-edge coloring}, Discrete Math. 341(3) (2018) 820--824.

\bibitem{WongZhu23Choos}
T. Wong, X. Zhu, \emph{Every graph is (2,3)-choosable}, Combinatorica  36(1) (2016) 121--127.

\bibitem{WongZhuChoos}
 T. Wong, X. Zhu, \emph{Total weight choosability of graphs}, J.
Graph Theory 66 (2011) 198--212.


\end{thebibliography}
\end{document}